\documentclass[12pt]{article}
\usepackage[top=0.8in,bottom=1in,left=0.75in,right=0.75in]{geometry}
\usepackage{amsmath,amssymb,amsfonts,epsf,euscript}
\usepackage{graphicx}
\usepackage[colorlinks=true]{hyperref}
\usepackage{float}
\newtheorem{theorem}{Theorem}[section]
\newtheorem{proposition}[theorem]{Proposition}
\newtheorem{lemma}[theorem]{Lemma}
\newtheorem{corollary}[theorem]{Corollary}

\newtheorem{definition}{Definition}[section]
\newtheorem{remark}{Remark}[section]
\newtheorem{proof}{P}
\newtheorem{Proof}{B}
\newtheorem{condition}{Condition}
\renewcommand{\theequation}{\arabic{section}.\arabic{subsection}.\arabic{equation}}
\title{\bf \Large Certain Semi-L\'evy Driven CARMA Processes: Estimation and Forecasting}
\author
{N. Modarresi\thanks{\scriptsize{
Department of Mathematics, Allameh Tabataba'i University, Tehran, Iran.
E-mail:  n.modarresi@atu.ac.ir.}}
\and  S. Rezakhah\thanks{\scriptsize{
Faculty of Mathematics and Computer Science, Amirkabir University of Technology, Tehran, Iran. Email: rezakhah@aut.ac.ir (S.Rezakhah) and m.mohammadiche@aut.ac.ir (M. Mohammadi).}}\and M. Mohammadi$^{\dagger}$}

\date{}
\begin{document}
\maketitle

\begin{abstract}
\noindent
Continuous-time autoregressive moving average (CARMA) process driven by simple semi-L\'evy process has periodically correlated property with many potential application in finance. In this paper, we study on the estimation of the parameters of the simple semi-L\'evy CARMA (SSLCARMA) process based on the Kalman recursion technique.
We implement this method in conjunction with the state-space representation of the associated process. The accuracy of estimation procedure is assessed in a simulated study. We fit a SSLCARMA(2,1) process to intraday realized volatility of Dow Jones Industrial Average data. Finally, We show that this process provides better in-sample forecasts of these data than the L\'evy driven CARMA process after deseasonalized them.\\


\end{abstract}

\section{Introduction}
Modeling of the continuous-time processes has a long history and has been carried out widely in financial econometrics. Early papers of Doop \cite{doop}, Phillips \cite{Phillips} and Durbin \cite{durbin} are dealt with properties and statistical analysis of Gaussian continuous-time ARMA (CARMA) processes. Brockwell \cite{bro1-2} introduced the L\'evy driven CARMA process for irregularly spaced data. These processes which are driven by non-decreasing L\'evy processes constitute a general class of stationary processes \cite{bro6}. Properties of second order L\'evy driven CARMA processes and some of their financial applications in modeling stochastic volatilities are discussed by Brockwell \cite{bro2}.
Strongly consistent estimators for the parameters of the subordinator CARMA processes based on uniformly spaced observations are presented in \cite{bro5}.\\
The L\'evy driven CARMA processes have the restriction that the underling process has stationary increments. In contrast, processes with periodically stationary increments such as semi-L\'evy processes have a wider application and are more prominent.
The semi-L\'evy processes have been extensively studied by Maejima and sato \cite{maj}. A class of CARMA processes driven by simple semi-L\'evy process which is denoted as SSLCARMA process, introduced by Modarresi et al. \cite{mod}. They studied the properties of this process and show that it is periodically correlated (PC).\\
In this paper, we study a certain class of CARMA($p, q$) process driven by simple semi-L\'evy compound Poisson process. In order to estimate parameters of the SSLCARMA($p, q$), first we characterize the sampled process. It is shown that the sampled process is a class of weak ARMA($p, p-1$) with independent and periodically identically distributed (ipid) noise. By the state-space representation of the sampled process, we compute the one-step linear prediction using Kalman recursion that is described in \cite{bbook2}. Numerical minimization of the sum of squares of errors gives least squares estimates of the parameters. The accuracy of estimation procedure is illustrated with simulated examples of some SSLCARMA(2,1) processes.\\
A growing number of research studies follow the intraday return that is determined by the availability of high-frequency financial data. Many of this data shows a PC structure in their squared log intraday returns \cite{r2}. For analysing such data, one approach is to remove the PC structure, then fit the corresponding stationary time series by the stationary process \cite{brodin}. The proposed SSLCARMA(2,1) process provides much better fitting to the 30-minute realized volatility series of 5-minute Dow Jones Industrial Average (DJIA) data which is applied by Brodin and Kl\"{u}ppelberg \cite{brodin}. For details on the determination of the realized volatility, see \cite{a}. We show the competitive performance of the SSLCARMA process with the L\'evy driven CARMA process. For this, we remove the periodicity of the 30-minute realized volatility series using filtering method \cite{brodin}, then fit a L\'evy driven CARMA process. Then we show that the SSLCARMA process forecast the sample paths of the 30-minute realized volatility much better than CARMA process.\\
The rest of the paper is organized as follows. In section 2, the definition and properties of the second order SSLCARMA process are reviewed.
We provide a discrete characterization of the SSLCARMA model through some proper discretization in section 3. The estimation of the parameters is followed by using the Kalman recursion algorithm to present one step ahead predictor model in this section as well. We show the performance of the estimation method by simulated data and also by applying the model to a real data set in section 4. Finally, we analysis the performance of the introduced model in compare with the L\'evy driven CARMA in some real data set. All proofs are given in Section 5.

\renewcommand{\theequation}{\arabic{section}.\arabic{equation}}
\section{Semi-L\'evy driven CARMA process}
\setcounter{equation}{0}
In order to define the simple semi-L\'evy driven CARMA, denoted by SSLCARMA process, first we present the simple semi-L\'evy (SSL) process. We remind that a semi-L\'evy process with period $T$ is a subclass of additive process with periodically stationary increments. Let $B_1, B_2, \ldots$ be a partition of the positive real line where $B_i=(s_{i-1}, s_i]$, $i\in\Bbb N$, $s_0=0$ and for some $r\in\Bbb N$, $T=\sum_{i=1}^{r}|B_i|$ where $|B_i|$ denotes the length of interval $B_i$. Also, $|B_i|=|B_{i+kr}|$ for $i,k\in\Bbb N$.
\begin{definition}\label{def 2.1}
The random measure $\{M(A): A\in{\cal B}\}$ where $\cal B$ is the Borel field on positive real line, is called simple semi-L\'evy (SSL) random measure with partition $B_1, B_2, \ldots$ and period $T=\sum_{i=1}^{r}|B_i|$ for some fixed $r\in\Bbb N$, if
\vspace{-1mm}
\begin{align*}
M\big(0,(k-1)T+s\big]=\sum_{i=1}^{(k-1)r+j-1}M_{_{i}}(s_{i-1},s_i]+M_{_{(k-1)r+j}}(s_{j-1},s],
\end{align*}
where $\{M_i, i\in\Bbb N\}$ is a sequence of L\'evy random measures that $M_{(k-1)r+i}$ is a copy of $M_i$ for all $i,k\in\Bbb N$. Moreover, $N(t):=M(0,t]$ is called SSL process.\\
If $\{M_i, i\in\Bbb N\}$ is a sequence of Poisson measures with rates $\lambda_i$ where $\lambda_{i+kr}=\lambda_i$, then $N(t):=M(0,t]$ is a SSL Poisson process with period $T$ and rate
\vspace{-1mm}
\begin{align}\label{Lam}
\Lambda_t=(k-1)\sum_{i=1}^{r}\lambda_i+\sum_{i=1}^{j-1}\lambda_i+ \frac{\lambda_j(s-s_{j-1})}{|B_{j}|},
\end{align}
for $t=(k-1)T+s$, $s\in B_j$. Therefore,
\vspace{-2mm}
\begin{align}\label{LI}
S_t=\gamma t+ \sum_{k=1}^{N(t)}J_k,
\end{align}
where $\gamma\in\Bbb R$ and $\{J_k: k\in\mathbb{N}\}$ is an independent and identically distributed (iid) sequence of random variables with probability distribution $F$ is called SSL compound Poisson process with drift.
\end{definition}
So, $E(S_t)=\gamma t+\Lambda_t\kappa$ and $\mbox{var}(S_t)=\Lambda_t\beta$ where $E(J_{k})=\kappa$, $E(J_{k}^{2})=\beta$ and $E(N(t))=\mbox{var}(N(t))=\Lambda_t$.
\begin{definition}\label{Definition2.2}
Let $\{S_t: t\geq0\}$ be a second order semi-L\'evy process with period $T$ defined by (\ref{LI}). The SSLCARMA$(p,q)$ process $\{Y_t:  t\geq0\}$, $p>q$, with parameters $a_1, \ldots, a_p, b_0, \ldots, b_q$ is the solution of the $p$th order stochastic differential equation $a(D)Y_t=b(D)DS_t$, where $D$ denotes differentiation with respect to $t$. The polynomials $a(z)=z^p+a_1z^{p-1}+ \ldots+a_p$ and $b(z)=b_0+b_1z+ \ldots+ b_{p-1}z^{p-1}$ have no common factors and the coefficients $b_j=0$ for $q<j<p$, $b_q=1$. The corresponding observation and state equations can be written as
\begin{align}\label{2}
Y_t=\bf b'X_t,
\end{align}
\vspace{-7mm}
\begin{align}\label{3}
d{\bf X}_t-A{\bf X}_tdt={\bf e}dS_t,
\end{align}
where
\begin{align*}
A=\begin{bmatrix}0&1&0&\ldots&0\\
0&0&1&\ldots&0\\
\vdots&\vdots&\vdots&\ddots&\vdots\\
0&0&0&\ldots&1\\
-a_p&-a_{p-1}&-a_{p-2}&\ldots&-a_1\\ \end{bmatrix},
&\hspace{1cm}{\bf e}=\begin{bmatrix}
0\\
0\\
\vdots\\
0\\
1\\
\end{bmatrix},
&\hspace{3mm}{\bf b}=\begin{bmatrix}
b_0\\
b_1\\
\vdots\\
b_{p-2}\\
b_{p-1}\\ \end{bmatrix}.
\end{align*}
\end{definition}
Every solution of equation $(\ref{3})$ satisfies the following relations for all $t>s\geq0$,
\begin{align}\label{4}
{\bf X}_t=e^{A(t-s)}{\bf X}_s+\int_s^te^{A(t-u)}{\bf e}dS_u,
\end{align}
where the paths of $S_{t}$ have bounded variation on compact intervals. From equation $(\ref{4})$ and the independence of the increments of $S_{t}$ one can easily verify that ${\bf X}_t$ is Markov. We characterized the moving average representation of the solution in (\ref{4}) and presented some properties of it. Furthermore, we show that the SSLCARMA process is verified to satisfy in some properties, if the following condition hold,
\begin{condition}\label{Condition 1}
The eigenvalues of the matrix $A$ have negative real parts and are distinct and the zeroes of the polynomial $a(z)$ are distinct. The assumption of distinct zeroes is not critical since multiple zeroes of $a(z)$ can be handled by replacing them with close but distinct zeroes and allowing each of these to converge to the multiple zero.
\end{condition}
In the following we extend the state process $\{\mathbf{X}_{t}: t\geq0\}$ to a process with index set $\mathbb{R}$. For this, we define the semi-L\'evy process on the whole real line.
\begin{definition}
Let $S_{t}$ be the semi-L\'evy process defined by (\ref{LI}). Then, the extend of the $S_{t}$ on the real line is defined as
\begin{align}\label{sl}
S_t:=S^{(1)}_{t}I_{[0,\infty)}(t)-S^{(2)}_{-t}I_{(-\infty,0]}(t),\qquad t\in\mathbb{R},
\end{align}
where $S_{t}^{(1)}$ and $S_{t}^{(2)}$ are independent copies of $S_{t}$, see \cite{mod}.
\end{definition}
\begin{remark}
Under Condition \ref{Condition 1} and $\lim_{t\rightarrow\infty}e^{At}=0$, as $s\rightarrow -\infty$, the solution (\ref{4}) with the specified properties satisfies
\vspace{-3mm}
\begin{align}\label{X}
{\bf X}_t=\int_{-\infty}^{t}e^{A(t-u)}{\bf e}dS_u.
\end{align}
\end{remark}

\begin{proposition}
If ${\bf X}_t$, defined by (\ref{X}), is independent of $\{S_r-S_t: r>t\}$ and the Condition \ref{Condition 1} holds, then the expected value and covariance function of ${\bf X}_t$ are periodic with period $T$ and consequently $Y_t$ is a periodically correlated (PC) process with period $T$, for $t\in\mathbb{R}$.
\end{proposition}
\vspace{-3mm}
For a proof, see \cite{mod}.
\begin{remark}{\ \\}
(i) If $S_t$ is a SSL process defined by (\ref{sl}) and the Condition \ref{Condition 1} is hold, then the SSLCARMA($p, q$) process with equations  (\ref{2}) and (\ref{X}) is defined as $Y_t=\mathbf{b}'\mathbf{X}_t=\int_{-\infty}^{\infty}\mathbf{b}'h(t-u)dS_{u},$
where $h(t)=e^{At}\mathbf{e}I_{[0,\infty)}(t)$ is called the kernel of the SSLCARMA process $Y_t$.\\ \\
(ii) If the kernel $h(\cdot)$ is non-negative and the jumps are additionally non-negative, then the process $(Y_{t})_{t\geq0}$ will be non-negative. The kernel is non-negative
if and only if the ratio $\frac{b(\cdot)}{a(\cdot)}$ is completely monotone \cite{tsai}, where the polynomials $a(\cdot)$ and $b(\cdot)$ is defined in Definition \ref{Definition2.2}. For SSLCARMA(2,1) process the condition is equivalent to the statement that the roots of $a(z)=0$, denoted by $\eta_{1}$ and $\eta_{2}$, are both real and that $b_{0}\geq min\{|\eta_{1}|, |\eta_{2}|\}$, \cite{bro5}.
\end{remark}

\section{Estimation procedure}
\setcounter{equation}{0}
In this section, we concerned with inference for the non-negative SSLCARMA process and deal with the problem of estimation of the parameters of this process. The theoretical properties of the corresponding time varying discrete-time process with equally spaced observation are developed. We apply an estimation method to estimate the coefficients of such sampled process which leads to estimate the parameters of SSLCARMA process.

\subsection{Characterization of the sampled process}
Following the method of Brockwell et al. \cite{bro5} and in order to estimate the parameters, we consider a discretization of the process. Let $\{Y_t: t\in\Bbb{R}^+\}$ be the SSLCARMA process with period $T>0$. We assume some equally spaced samples as $\{Y_n:=Y_{nh}, n=1,\ldots, N\}$ where $h=\frac{T}{M_0}$ and $N\in\mathbb{N}$. It is shown in \cite{mod} that $\{Y_t, t\in\Bbb{R}\}$ is a PC with period $T$, so $\{Y_n, n=1,\ldots, N\}$ is PC with period $M_0$. Therefore, we have the following result.

\begin{proposition}
Let $\{Y_{n}: n=1,\ldots, N\}$ be the available sampled SSLCARMA$(p,q)$ process with period $M_0$, then under Conditions \ref{Condition 1} the process can be written as $Y_n=\sum_{r=1}^{p}Y_{n}^{(r)}$ where
\vspace{-1mm}
\begin{align}\label{u2}
Y_{n}^{(r)}&=\int_{-\infty}^{nh}\alpha_re^{\eta_r(nh-u)}dS_u,
\end{align}
$\alpha_r=\frac{b(\eta_r)}{a'(\eta_r)}$ in which $\eta_{r}$ and $a'(\cdot)$ are the roots and the derivative of the autoregressive polynomial $a(\cdot)$ presented in Definition \ref{Definition2.2}, respectively.
\end{proposition}
The proof, which is an immediate result of the decomposition of the integrand in (\ref{X}) is the same as the one presented \cite{bro5}.

\begin{corollary}
For $n\in\Bbb N$, a closed formula for the sampled process $Y_{n}^{(r)}$ is
\begin{align}\label{u3}
Y_{n}^{(r)}=e^{\eta_rh}Y_{n-1}^{(r)}+Z_n^{(r)},
\end{align}
where $Z_n^{(r)}=\alpha_r\int_{(n-1)h}^{nh}e^{\eta_r(nh-u)}dS_u$ is an independent and periodically identically distributed (ipid) noise.
\end{corollary}
\vspace{-2mm}
Proof: see Appendix A, P\ref{P1}.\\ \\
Now in the following lemma, assuming some conditions on autocovariance function, we show that any PC process can be represented as a moving average process with ipid noise. So, this lemma can be applied to the sampled SSLCARMA process which has been proved in \cite{mod} that is PC process and leads to a class of weak ARMA process with ipid noise.

\begin{lemma}\label{MA}
Let $\{G_n: n\in\Bbb N\}$ be a zero-mean PC process with period $M_0\in\Bbb N$ and $\gamma_n(l)=cov(G_n, G_{n+l})=0$ while $l$ is greater than some integer $p$. Then $\{G_n: n\in\Bbb N\}$ can be represented as a moving average process with ipid noise of order $p$ with constant coefficients of some uncorrelated and PC random variables $\{\xi_n: n\in\Bbb N\}$ as
\begin{align*}
G_n=\xi_n+\theta_1\xi_{n-1}+\ldots+\theta_p\xi_{n-p}.
\end{align*}
\end{lemma}
\vspace{-1mm}
Proof: see Appendix A, P\ref{P2}.

\begin{theorem}
Let $\phi(B)=\prod_{i=1}^{p}(1-e^{\eta_ih}B)=: 1-\varphi_1B-\varphi_2B^2-\ldots-\varphi_pB^p$ be an operator and $B^jY_n=Y_{n-j}$. By applying $\phi(B)$ to each elements of $Y_n=\sum_{r=1}^{p}Y_{n}^{(r)}$ and summing over $r$, the sampled SSLCARMA process $Y_n$ yields to
\begin{align}\label{u6}
\phi(B)Y_n=V_{n}^{(1)}+V_{n-1}^{(2)}+\ldots+V_{n-p+1}^{(p)},
\end{align}
where for each fixed $k\in\{1, \ldots, p\}$, $V_{n-k+1}^{(k)}$ is an ipid sequence with period $M_0$ defined by
\begin{align*}
V_{n-k+1}^{(k)}=\sum_{r=1}^{p}\big(e^{(k-1)\eta_{r}h}-\sum_{j=1}^{k-1}\varphi_{j}e^{(k-1-j)\eta_{r}h}\big)\alpha_r\int_{(n-k)h}^{(n-k+1)h}
e^{\eta_r((n-k+1)h-u)}dS(u).
\end{align*}
\end{theorem}
\vspace{-1mm}
Proof: see Appendix A, P\ref{P3}.

\begin{remark}\label{remark3.1}
It follows from (\ref{u6}) that $\phi(B)Y_{n}=V_{n}^{(1)}+V_{n-1}^{(2)}+\ldots+V_{n-p+1}^{(p)}$ is a $(p-1)-$dependent sequence. So, if the PC process $Y_n$ has zero mean, then by Lemma \ref{MA} there exists an uncorrelated PC noise $\xi_n$ such that $\phi(B)Y_n$ is moving average process with ipid noise of order $(p-1)$ in which
\begin{align*}
\phi(B)Y_n&=\xi_n+\sum_{j=1}^{p-1}\theta_j\xi_{n-j}=(1+\sum_{j=1}^{p-1}\theta_jB^j)\xi_n
=(\sum_{j=0}^{p-1}\theta_jB^j)\xi_n=:\theta(B)\xi_n,
\end{align*}
where $\theta_{0}:=1$ and coefficients $\theta_j$, $j=1, \ldots, p-1$, are constant depending on the parameters of the SSLCARMA process. Therefore
\begin{equation}\label{u8}
\phi(B)Y_n=\theta(B)\xi_n.
\end{equation}
\end{remark}
So, $Y_n$ is a class of weak ARMA$(p,p-1)$ process with ipid noise and from (\ref{u6}), $\xi_n=\theta(B)^{-1}\sum_{r=1}^{p}V_{n-r+1}^{(r)}.$

\subsection{Kalman prediction}
The Kalman filter is an optimal estimating method that infers parameters from indirect and uncertain observations. It is recursive so that new measurements can be processed as they arrive. This method minimizes the mean square error of the estimated parameters. For more details see \cite{bbook2}, chapter 9. In order to present a prescription of the optimal filter we find the prerequisites of the algorithm such as the covariance matrix of the noise and linear predictors. \\
By Remark \ref{remark3.1}, the centered sampled process $Y^{*}_n:=Y_n-E(Y_n)$, where $E(Y_{n})$ is period mean, satisfies the class of weak ARMA process (\ref{u8}) driven by ipid noise $\xi_n^{*}=\xi_n-E(\xi_n)$. It follows also from (\ref{2}) that the process $Y^{*}_n$ has the observation equation
\begin{align}\label{FE1}
Y^{*}_n={\bf b}'{\bf X}^{*}_n,
\end{align}
where ${\bf X}^{*}_n:={\bf X}_n-E({\bf X}_n)$ is the centered state vector of ${\bf X}_n:={\bf X}_{nh}$. It satisfies the state equation
\begin{align}\label{FE2}
{\bf X}^{*}_n=e^{Ah} {\bf X}^{*}_{n-1}+{\bf U}_n,
\end{align}
where ${\bf U}_{n}:=\int_{(n-1)h}^{nh}e^{A(nh-u)}{\bf e}dS_u+e^{Ah}E({\bf X}_{n-1})-E({\bf X}_n)$ is a sequence of zero-mean ipid random vectors with covariance matrices
\begin{align}\label{mat}
Q_{n}=
\begin{cases}
\frac{\beta\lambda_j}{|B_{j}|}\int_{0}^{h}e^{Au}{\bf e}{\bf e'}e^{A'u}du, &  (n-1)h, nh\in B_{j}
\\
\\
\frac{\beta\lambda_j}{|B_{j}|}\int_{nh-s_{j}}^{h}e^{Au}{\bf e}{\bf e'}e^{A'u}du
+\frac{\beta\lambda_{j+1}}{|B_{j+1}|}\int_{0}^{nh-s_{j}}e^{Au}{\bf e}{\bf e'}e^{A'u}du, & (n-1)h\in B_{j}, nh\in B_{j+1}
\end{cases}
\end{align}
in which $j\in\mathbb{N}$, $\beta=E(J_{k}^{2})$ and $\lambda_{1}, \ldots, \lambda_{r}$ are jump-rates corresponding to the increments of the SSL Poisson random measure on partitions $B_{1}, \ldots, B_{r}$, as assumed in Definition \ref{def 2.1}. For more details regarding (\ref{FE2}) and (\ref{mat}), see  Appendix B, B\ref{B1}.
\begin{remark}\label{remark3.2}
Since the sampled process $\{Y_{n}: n=1, \ldots, N\}$ is a PC process with period $M_0$, the periodic mean $E(Y_{n})$ is estimated by sample periodic mean
\begin{align*}
\overline{Y}_n=\frac{1}{[\frac{N}{M_{0}}]}\sum_{i=1}^{[\frac{N}{M_{0}}]}Y_{m+(i-1)M_0},
\end{align*}
where $m=n-[\frac{n}{M_{0}}]M_{0}$ in which $[x]$ denotes the integer part of $x$, see \cite{hurd}, Chapter 9.
\end{remark}
The inferential goal is to estimate of the SSLCARMA parameter vector $(a_{1}, \ldots, a_{p}, b_{0},$ $\ldots, b_{q-1})^{\prime}$. We do this by using the Kalman recursions in conjunction with the state-space representation in equations (\ref{FE1}) and  (\ref{FE2}). We compute the one-step linear predictors $\widehat{Y}_{n}^{*}:=P_{n-1}(Y_n^*)$ in terms of $Y_0^*, \ldots, Y_{n-1}^*$, $n=1, \ldots, N$, based on the Kalman algorithm which is summarized in Table \ref{table1}. Numerical minimization of the sum of squares of these one-step errors, $\sum_{n=1}^{N}\big(Y_n^*-\widehat{Y}_{n}^{*}\big)^{2}$, with respect to parameters of the model gives least squares estimates of the SSLCARMA coefficients. The Kalman filter algorithm can be roughly organized under the following steps.
\vspace{-3mm}
\begin{table}[h]
\caption{Kalman Recursion Algorithm}
\vspace{2mm}
\centering
\begin{tabular}{c}
\hline
\begin{minipage}{6.7in}
\vspace{1mm}
\begin{itemize}
\item[(a)] The predictors $\widehat{Y}_{n}^{*}$ of the state-space model (\ref{FE1}) and (\ref{FE2}) are determined by the one-step predictors $\widehat{\mathbf{X}}_{n}^{*}:=P_{n-1}\mathbf{X}_{n}^{*}$, the error covariance matrices $\Omega_{n}:=E\big[(\mathbf{X}_{n}^{*}-\widehat{\mathbf{X}}_{n}^{*})(\mathbf{X}_{n}^{*}-\widehat{\mathbf{X}}_{n}^{*})'\big]$ and the initial conditions
\begin{itemize}
\item[(i)] $\widehat{\mathbf{X}}_{1}^{*}={\bf 0}\Longrightarrow\widehat{Y}_{1}^{*}=\mathbf{b}'\widehat{\mathbf{X}}_{1}^{*}=0$
\item[(ii)] $\Omega_{1}=\frac{\beta\lambda_{1}}{|B_{1}|}\int_{0}^{h}e^{Au}\mathbf{e}
\mathbf{e}^{\prime}e^{A^{\prime}u}du+\sum_{j=1}^{r}
\sum_{k=0}^{\infty}\frac{\beta\lambda_{j}}{|B_{j}|}\int_{s_{_{j-1}}+h+kT}^{s_{_{j}}+h+kT}e^{Au}\mathbf{e}\mathbf{e}^{\prime}e^{A^{\prime}u}du$,
\end{itemize}
where $T=hM_0$. For more details regarding $\Omega_{1}$, see Appendix B, B\ref{B2}.
\item[(b)] For $n=1, \ldots, N$
\begin{itemize}
\item[(i)] $\widehat{\mathbf{X}}_{n+1}^{*}=e^{Ah}\widehat{\mathbf{X}}_{n}^{*}+\Theta_{n}\Delta_{n}^{-1}\big(Y_{n}^{*}-\widehat{Y}_{n}^{*}\big) \Longrightarrow \widehat{Y}_{n+1}^{*}=\mathbf{b}'\widehat{\mathbf{X}}_{n+1}^{*}$,
\item[(ii)] $\Omega_{n+1}=e^{Ah}\Omega_{n}e^{A'h}+Q_{n}-\Theta_{n}\Delta_{n}^{-1}\Theta'_{n}$
\end{itemize}
\noindent where $Q_{n}$ is defined in (\ref{mat}), $\Delta_{n}=\mathbf{b}'\Omega_{n}\mathbf{b}$, $\Theta_{n}=e^{Ah}\Omega_{n}\mathbf{b}$
and $\Theta'_{n}$ is the transpose of the vector $\Theta_{n}$.
\end{itemize}
\end{minipage}
\vspace{2mm}
\\
\hline
\end{tabular}
\label{table1}
\end{table}

\begin{remark}
In this algorithm we assume that $\gamma=0$ and $\kappa$ and $\beta$ are not represented in the predictor $\widehat{\mathbf{X}}_{n}^{*}$ but also in $\widehat{Y}_{n}^{*}$. Moreover, the parameter $\beta$ in $\widehat{\mathbf{X}}_{n}^{*}$ is omitted. Then by minimizing the sum of the square errors with respect to $(a_{1}, \ldots, a_{p}, b_{0},\ldots, b_{q-1}, \lambda_{1}, \ldots, \lambda_{r})^{\prime}$ the coefficients parameters of SSLCARMA process and $\lambda_i$, $i=1, \ldots, r$ are estimated.
\end{remark}

\section{Data analysis}
In this section, we conduct a simulation study to test the estimation procedure for the SSLCARMA parameters and to assess the quality of the estimates in subsection 4.1. Then in subsection 4.2, we apply the introduced model to the intraday realized volatility data for the Dow Jones Industrial Average (DJIA) and by minimization the sum of squared errors we estimate the parameters of the process. The estimation results are compared by the model which is introduced by Brodin and Kl\"{u}ppelberg \cite{brodin} that is fitted by a L\'evy-driven CARMA process after removing its periodicity in subsection 4.3.

\subsection{Simulation study}
\begin{figure}[H]
\centering
\includegraphics[scale=0.53]{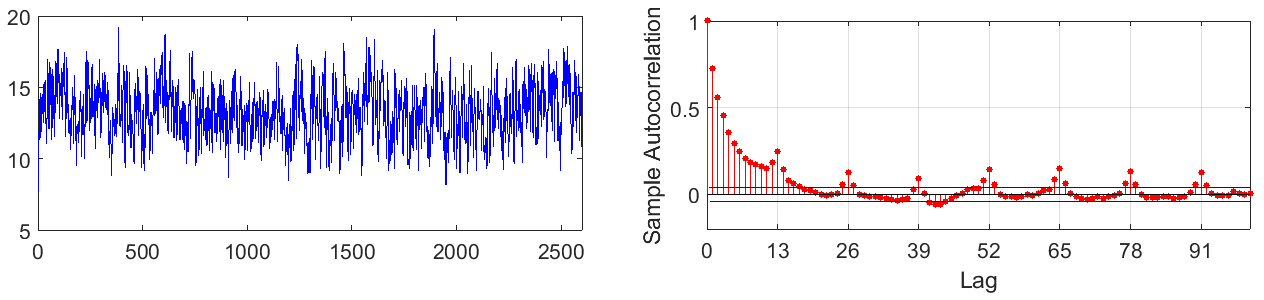}
\vspace{-3mm}
\caption{\scriptsize{The sample paths of the process (left) and the sample ACF (right) of the SSLCARMA(2,1).}}
\label{figure1}
\end{figure}

\setcounter{equation}{0}
The simulation of SSLCARMA process is followed by the moving average representation (\ref{X}). For this, first we simulate the SSL process $S_t$ represented in (\ref{LI}). As a special case, we assume $S_t$ to be a SSL process with period $T=13$ and the lengths of the successive subintervals of each period interval are considered as $10, 2$ and $1$. The arrival rates of the semi-L\'evy Poisson process on these subintervals are assumed as $\lambda_{1}=10, \lambda_{2}=15$ and $\lambda_{3}=3.$  Moreover, the jumps $J_{k}$ are assumed to be exponentially distributed with parameter $\eta=0.25$. We simulate 1000 realizations of the SSLCARMA(2,1) that is specified by the equation
$$(D^2+a_1D+a_2)Y_t=(b_0+D)DS_t,$$
where $b(z)=b_0+z$, $a(z)=z^2+a_1z+a_2$ and $S_t$ is the SSL compound Poisson process with drift. In this case, the parameters are $a_{1}=3, a_{2}=0.5$ and $b_{0}=2$ and the simulation is for the duration of 200 period intervals. Then, each realization is sampled at spaces $h=1$. Figure \ref{figure1}, shows the sample path and sample autocorrelation function (ACF) of this SSLCARMA(2,1) process.
\vspace{-1mm}
\noindent
For each realization we compute least squared estimators of the parameters of the SSLCARMA(2,1) process.
As noted in Remark 3.3, it is not required to estimate the parameter of jump distribution $\eta$. The sample mean, bias and standard deviation of these estimators are shown in Table \ref{table2}.
\begin{table}[h]
\caption{Simulation study based on 1000 realizations for the SSLCARMA(2,1) parameters.}
\centering
\begin{tabular}{c c c c c c c c c c c c c c c c c c}
\hline\hline
& & $a_{1}$ & & & $a_{2}$ & & & $b_{0}$ & & & $\lambda_{1}$ & & & $\lambda_{2}$ & & & $\lambda_{3}$\\
True & & 3 & & & 0.5 & & & 2 & & & 10 & & & 15 & & & 3\\
\hline
Mean & & 2.9146 & & & 0.5026 & & & 2.0106  & & & 10.1118 & & & 14.8977 & & & 3.1114 \\
Bais & & 0.0854 & & & 0.0026 & & & 0.0106 & & & 0.1118 & & & 0.1023 & & & 0.1114 \\
Std. dev. & & 0.0611 & & & 0.0145 & & & 0.0486 & & & 0.0269 & & & 0.0250 & & & 0.0222 \\
\hline
\end{tabular}
\label{table2}
\end{table}

\subsection{Intraday realized volatility for the DJIA}
\begin{figure}[H]
\centering
\includegraphics[scale=0.53]{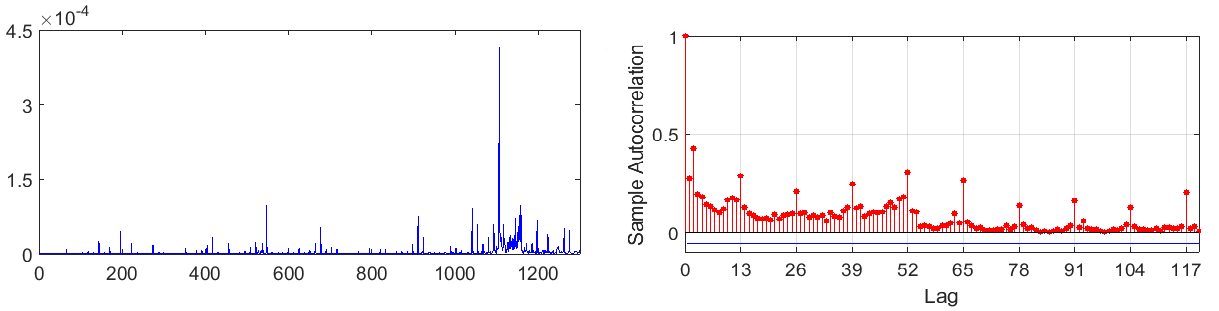}
\vspace{-3mm}
\caption{\scriptsize{The sample paths of the process (left) and the sample ACF (right) of $\{RV_{n}: 1, \ldots, 1300\}$.}}
\label{figure2}
\end{figure}

\begin{figure}[H]
\centering
\includegraphics[scale=0.5]{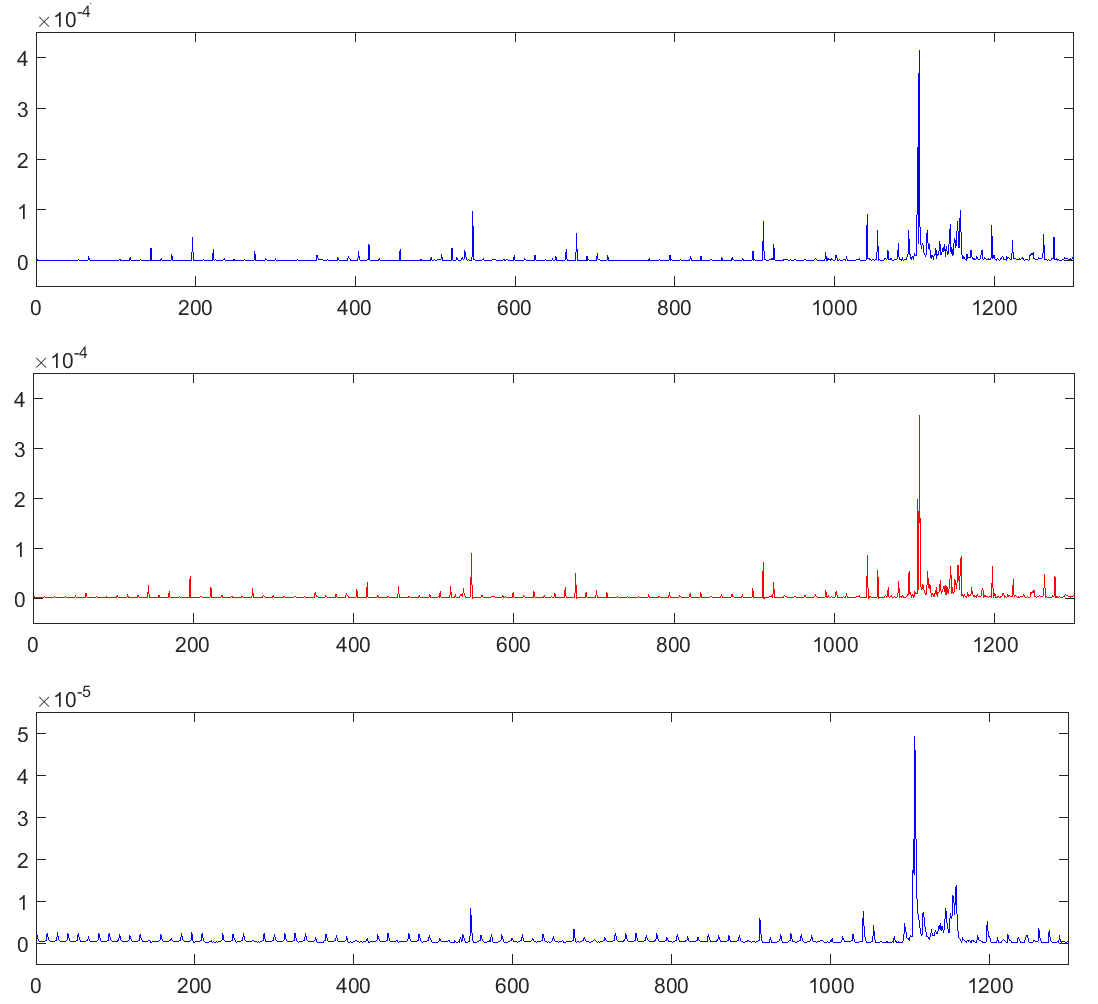}
\vspace{-3mm}
\caption{\scriptsize{The time series $RV_{n}$ (top), the $\widehat{RV}_{n}$ (middle) and the absolute errors of the predictor $\widehat{RV}_{n}$ (bottom).}}
\label{figure3}
\end{figure}

\setcounter{equation}{0}
Realized volatility is a non-parametric estimate of the return variation. The most obvious realized volatility measure is the sum of finely-sampled squared return realizations over a fixed time interval as
\vspace{-4mm}
\begin{align}\label{RV}
RV_{n}=\sum_{j=1}^{k}d_{n,j}^{^{2}},
\end{align}
where $d_{n,j}=ln(P_{_{j+1+k(n-1)}})-ln(P_{_{j+k(n-1)}})$ in which $P_{_{j}}$ is asset price. One of the striking features of financial time series is that the log returns have negligible correlation while its squared log returns are significantly correlated \cite{b5}. Many of high-frequency time series show a PC structure in their squared log returns \cite{r2}, so according to relation  (\ref{RV}) the intraday realized volatility have PC structure.\\
Here, we describe the application of the estimation procedure to the 30-minute realized volatility, denoted by $RV_{n}$, of the 5-minute DJIA data. This data set is recorded between 9:35 to 16:00 from October 3th, 2017 to February 27th, 2018. There was a total of $N=100$ trading days not including the weekends and holidays with 78 5-minute observations per day, resulting in the total of 7800 5-minute observations. We compute the $RV_{n}$ form these data by (\ref{RV}) for $k=6$. Figure \ref{figure2} shows the sample paths and the sample ACF of the time series $\{RV_{n}: n=1, \ldots, 1300\}$.
\vspace{-3mm}
\vspace{-1mm}
\noindent
It is clear from ACF that the time series $RV_{n}$ have a PC structure with period 13. We fit a SSLCARMA(2,1) model to $RV_{n}$. For this, we consider the SSL process $S_{t}$, defined by (\ref{sl}), as the underlying process with period $M_0=13$. Furthermore, the lengths of the successive subintervals of each period interval are 10, 2, 1 where corresponding arrival rates of the semi-L\'evy Poisson process on these subintervals are $\lambda_{1}, \lambda_{2}$ and $\lambda_{3}$ respectively. We use the Kalman algorithm which presented in Table \ref{table1} and compute the one-step predictions $\widehat{RV}^{*}_{n}:=P_{n-1}(RV_n^*)$.
By numerical minimization of the sum of squared errors, $\sum_{n=1}^{1300}(RV^{*}_{n}-\widehat{RV^{*}_{n}})^{2}$, where the centered realized volatility $RV^{*}_{n}=RV_{n}-\overline{RV_{n}}$ in which $\overline{RV_n}$ is followed from Remark \ref{remark3.2}, we estimate the parameters of the SSLCARMA(2,1) process. Table \ref{table3}, shows the outcomes of estimating the parameters of the SSLCARMA(2,1).
\vspace{-3mm}
\begin{table}[H]
\caption{Estimated parameters of the SSLCARMA(2,1).}
\centering
\begin{tabular}{c c c c c c c c c c c c c c c c}
\hline\hline
$\hat{a}_{1}$ & & $\hat{a}_{2}$ & & $\hat{b}_{0}$ & & $\hat{\lambda}_{1}$ & & $\hat{\lambda}_{2}$ & & $\hat{\lambda}_{3}$\\
\hline
$1.0472$ & & $0.2158$ & & $1.0843$ & & $3.5099$ & & $6.2535$ & & $14.3454$\\
\hline
\end{tabular}
\label{table3}
\end{table}
\vspace{-3mm}
\noindent
We use these estimators in state-space representation (\ref{FE1}) and (\ref{FE2}) and compute the one-step Kalman predictions of the  $\widehat{RV^{*}_{n}}$. So, we can be compute the one-step predictions $\widehat{RV_{n}}\approx\widehat{RV^{*}_{n}}+\overline{RV_{n}}$. Figure \ref{figure3} shows the time series $RV_{n}$ and the $\widehat{RV_{n}}$ for $n=1, \ldots, 1300$.
\vspace{-1mm}
\vspace{-5mm}
\subsection{In-sample performance analysis}
\begin{figure}[H]
\centering
\includegraphics[scale=0.5]{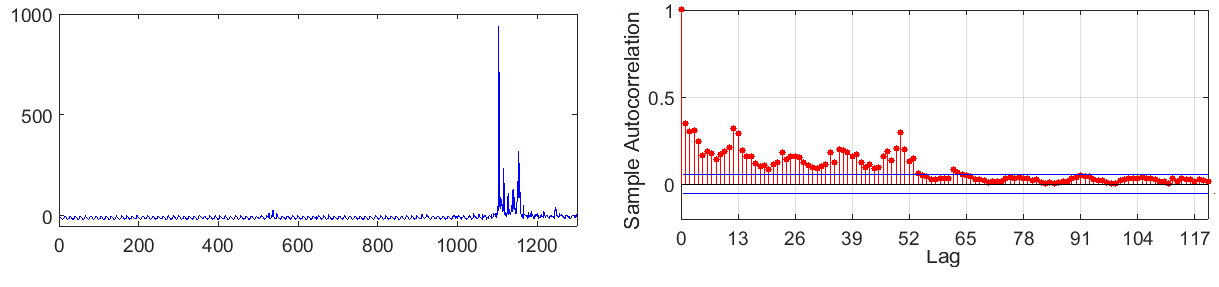}
\vspace{-3mm}
\caption{\scriptsize{The sample paths of the process (left) and the sample ACF (right) of $\{rv_{n}: n=1, \ldots, 1300\}$.}}
\label{figure4}
\end{figure}

\begin{figure}[H]
\centering
\includegraphics[scale=0.5]{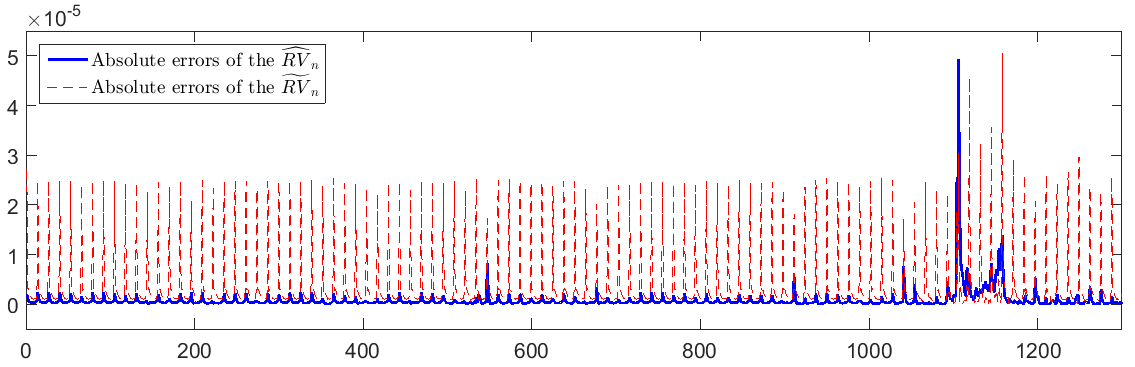}
\vspace{-3mm}
\caption{\scriptsize{The absolute errors of the $\widehat{RV}_{n}$ (blue) and of the $\widetilde{RV}_{n}$ (red).}}
\label{figure5}
\end{figure}

To compare the performance of the SSLCARMA process with a L\'evy driven CARMA process, we consider one-step Kalman prediction errors for the PC time series $RV_{n}$ applied in subsection 4.2.\\
For modeling $RV_{n}$ using the L\'evy driven CARMA process, we follow the method of Brodin and Kl\"{u}ppelberg \cite{brodin}, and remove the period of this time series using filtering method in \cite{brodin},
\begin{align}\label{rv}
rv_{n}:=\frac{RV_{n}-\hat{\mu}}{\hat{\nu}_{n}},\qquad n=1, \ldots, 1300,
\end{align}
where $\hat{\mu}$ is the sample mean of $RV_{n}$ and $\hat{\nu}_{n}$ are the seasonality coefficients estimated by
\begin{align*}
\hat{\nu}_{n}=median_{_{i=1, 2, \ldots, 100}}|RV_{m+(i-1)M_0}|,\qquad n=1, \ldots, 1300,
\end{align*}
in which $M_0=13$ and $m=n-[\frac{n}{M_0}]M_0$. Figure \ref{figure4} shows the sample paths and the sample ACF of the filtered time series $\{rv_{n}: n=1, \ldots, 1300\}$ which is defined by (\ref{rv}). As it
is shown, the filtered data has no clear periodicity effect.
\noindent
First, we model the time series $rv_{n}$ by a L\'evy driven CARMA($p,q$) process using the method of the paper \cite{bro5}, based on the assumption that the L\'evy process $\{L_{t}: t\in\mathbb{R}\}$ is a compound Poisson process with arrival rate $\lambda$ and exponentially distributed jump size.
It follows that the sampled process $Y_{n}:=Y_{nh}$ is the weak ARMA($p,p-1$) process driven by the white noise sequence. So, the centered sampled process $Y^{*}_n:=Y_n-E(Y_n)$ satisfies a weak ARMA process with state-space representation
\begin{align}\label{state2}
Y^{*}_n={\bf b}'{\bf X}^{*}_n \qquad \text{and} \qquad {\bf X}^{*}_n=e^{Ah} {\bf X}^{*}_{n-1}+{\bf U}_n,
\end{align}
where $\mathbf{U}_{n}=\int_{(n-1)h}^{nh}e^{A(nh-u)}{\bf e}dL_u+e^{Ah}E({\bf X}_{n-1})-E({\bf X}_n)$ is a sequence of zero-mean iid random vectors with covariance matrix
$Q_n=\text{var}(L_{1})\int_{0}^{h}e^{Au}\mathbf{e}\mathbf{e}'e^{A'u}du$ (for more details see Appendix B, B\ref{B3}).
We use the state-space equation (\ref{state2}) with error covariance matrix  $\Omega_{1}=\text{var}(L_{1})\int_{0}^{\infty}e^{Au}\mathbf{e}\mathbf{e}'e^{A'u}du$ (provided in Appendix B, B\ref{B4}) in Kalman algorithm which presented in Table \ref{table1} and compute the one-step predictions $\widetilde{rv}^{*}_{n}:=P_{n-1}(rv_n^*)$. By numerical minimization of the sum of squared errors, $\sum_{n=1}^{1300}(rv^{*}_{n}-\widetilde{rv}^{*}_{n})^2$ where $rv^{*}_{n}=rv_{n}-\overline{rv_{n}}$ in which $\overline{rv_{n}}=\frac{1}{1300}\sum_{i=1}^{1300}rv_{i}$, we estimate the parameters of the CARMA(2,1) process as $\hat{\alpha}_{1}=0.3292, \hat{\alpha}_{2}=0.0137$ and $\hat{\beta}_{0}=0.2250$. We use these estimators in state-space equation (\ref{state2}) and compute the corresponding
one-step Kalman predictions of the $\widetilde{rv}^{*}_{n}$. So, from this and the filtering method (\ref{rv}), we predict the intraday realized volatility of the main data as $\widetilde{RV}_{n}$ by the followings. So,
\begin{align*}
\widetilde{RV}_{n}\approx\hat{\nu}_{n}\widetilde{rv}_{n}
+\hat{\mu}, \qquad n=1, \ldots, 1300.
\end{align*}
In Figure \ref{figure5}, we illustrate the absolute errors of the one-step Kalman predictions $\widehat{RV}_{n}$, which is computed from SSLCARMA(2,1), and $\widetilde{RV}_{n}$. The mean absolute error of the $\widehat{RV}_{n}$ and $\widetilde{RV}_{n}$, respectively, are $8.1390\times10^{-7}$ and $3.4732\times10^{-6}$.

\section{Appendix}
{\bf\large Appendix A}
\setcounter{equation}{0}
\begin{proof}\label{P1}
\textnormal{\textbf{: Proof of Corollary 3.2}}
\end{proof}
\vspace{-3mm}
In discrete form, we can rewrite $Y_n^{(r)}$ in (\ref{u2}) as
\begin{align*}
Y_{n}^{(r)}&=\int_{-\infty}^{(n-1)h}\alpha_re^{\eta_r(nh-u)}dS_u+\int_{(n-1)h}^{nh}\alpha_re^{\eta_r(nh-u)}dS_u\nonumber\\
&=e^{\eta_r h}\int_{-\infty}^{(n-1)h}\alpha_re^{\eta_r((n-1)h-u)}dS_u+\int_{(n-1)h}^{nh}\alpha_re^{\eta_r(nh-u)}dS_u
=:e^{\eta_r h}Y_{n-1}^{(r)}+Z_n^{(r)}.
\end{align*}
\vspace{-2mm}
\begin{proof}\label{P2}
\textnormal{\textbf{: Proof of Lemma 3.3}}
\end{proof}
We define the subspace ${\cal{M}}_n=\overline{sp}\{G_m,-\infty<m\leq n\}$ of $L^2$ for each $n\in\Bbb N$ and set
\begin{equation}\label{ma2}
\xi_n=G_n-P_{{\cal{M}}_{n-1}}G_n,
\end{equation}
where $P_{{\cal{M}}_{n-1}}$ denote the projection mapping onto subspace ${\cal{M}}_{n-1}$. Clearly $\xi_n\in{\cal{M}}_n$ and $\xi_n\in{\cal{M}}_{n-1}^\bot$ where ${\cal{M}}_{n-1}^\bot$ is orthogonal complement of subset ${\cal{M}}_{n-1}$. Thus for $m<n$, $\xi_m\in{\cal{M}}_m\subset{\cal{M}}_{n-1}$ and $\xi_n\in{\cal{M}}_{n-1}^\bot$ and hence $E[\xi_n\xi_m]=0$. For more details regarding the space $L^{2}$ and projection, see \cite{bbook}, chapter 2. Furthermore, we can show that
\begin{align*}
\lim_{s\rightarrow\infty}P_{\overline{sp}\{G_m, m=n-s, \ldots, n-1\}}G_n=P_{{\cal{M}}_{n-1}}G_n.
\end{align*}
Since $\{G_n, n\in\Bbb N\}$ is PC with period $M_0$ and $L^2$ norm, $\|\xi_{n}\|=\sqrt{E(\xi_{n}^{2})}$, is continuous, we have that
\begin{align*}
\|\xi_{n+M_0}\|&=\|G_{n+M_0}-P_{{\cal{M}}_{n+M_0-1}}G_{n+M_0}\|=\lim_{s\rightarrow\infty}\|G_{n+M_0}-P_{\overline{sp}\{G_m, m=n+M_0-s, \ldots, n+M_0-1\}}G_{n+M_0}\|\\
&=\lim_{s\rightarrow\infty}\|G_{n}-P_{\overline{sp}\{G_m, m=n-s, \ldots, n-1\}}G_{n}\|=\|G_n-P_{{\cal{M}}_{n-1}}G_n\|=\|\xi_n\|.
\end{align*}
We conclude that $\xi_n$ is a zero-mean noise with variance $\sigma^2_n=\|\xi_n\|^2$ where $\sigma^2_{n+M_0}=\sigma^2_n$.
Now by (\ref{ma2}), it follows that
\begin{align*}
{\cal{M}}_{n-1}=\overline{sp}\{G_m, m<n-1, \xi_{n-1}\}=\overline{sp}\{G_m, m<n-p, \xi_{n-p}, \ldots, \xi_{n-1}\}.
\end{align*}
The subspace ${\cal{M}}_{n-1}$ can be decomposed into the two orthogonal subspaces $\overline{sp}\{\xi_{n-p}, \ldots, \xi_{n-1}\}$ and ${\cal{M}}_{n-p-1}$. Since by the assumption $\gamma_n(l)=0$ for $|l|>p$, therefore $G_n\bot{\cal{M}}_{n-p-1}$ and so by the properties of projection mappings and Theorem 2.4.1 in \cite{bbook} we have
\begin{align*}
P_{{\cal{M}}_{n-1}}G_n&=P_{{\cal{M}}_{n-p-1}}G_n+P_{\overline{sp}\{\xi_{n-p}, \ldots, \xi_{n-1}\}}G_{n}\\
&=0+\sigma^{-2}E[G_n\xi_{n-1}]\xi_{n-1}+ \ldots+ \sigma^{-2}E[G_n\xi_{n-p}]\xi_{n-p},
\end{align*}
and by denoting $\theta_j:=\sigma^{-2}E[G_n\xi_{n-j}]$ and substituting $P_{{\cal{M}}_{n-1}}G_n$ in (\ref{ma2}) we have
$G_n-\xi_n=\theta_1\xi_{n-1}+ \ldots+ \theta_p\xi_{n-p}$.

\begin{proof}\label{P3}
\textnormal{\textbf{: Proof of Theorem 3.4}}
\end{proof}
It follows from $Y_n=\sum_{r=1}^{p}Y_{n}^{(r)}$ that
\begin{align*}
\phi(B)Y_n&=\sum_{r=1}^{p}\phi(B)Y_{n}^{(r)}=\sum_{r=1}^{p}\prod_{i=1}^{p}(1-e^{\eta_ih}B)Y_{n}^{(r)}\\
&=\sum_{r=1}^{p}\prod_{i\neq r}(1-e^{\eta_ih}B)(1-e^{\eta_{r}h}B)Y_{n}^{(r)}=\sum_{r=1}^{p}\prod_{i\neq r}(1-e^{\eta_ih}B)(Y_{n}^{(r)}-e^{\eta_{r}h}Y_{n-1}^{(r)})
\end{align*}
By denoting $1+\psi_{1}B+\psi_{2}B^{2}+\ldots+\psi_{p-1}B^{p-1}:=\prod_{i\neq r}(1-e^{\eta_ih}B)$ and by (\ref{u3}), we have
\vspace{-2mm}
\begin{align}\label{The1}
\phi(B)Y_n=\sum_{r=1}^{p}\big(1+\psi_{1}B+\psi_{2}B^{2}+\ldots+\psi_{p-1}B^{p-1}\big)Z_{n}^{(r)}.
\end{align}
Since $(1-e^{\eta_{r}h}B)\prod_{i\neq r}(1-e^{\eta_ih}B)=\prod_{i=1}^{p}(1-e^{\eta_ih}B)=\phi(B)$, it follows from assumption of the theorem that
\begin{align*}
1-&\varphi_{1}B-\varphi_{2}B^{2}-\ldots-\varphi_{p}B^{p}=(1-e^{\eta_{r}h}B)\prod_{i\neq r}(1-e^{\eta_ih}B)\\
&=(1-e^{\eta_{r}h}B)(1+\psi_{1}B+\psi_{2}B^{2}+\ldots+\psi_{p-1}B^{p-1})\\
&=1+(\psi_{1}-e^{\eta_{r}h})B+(\psi_{2}-\psi_{1}e^{\eta_{r}h})B^{2}+\ldots+(\psi_{p-1}-\psi_{p-2}e^{\eta_{r}h})B^{p-1}-\psi_{p-1}e^{\eta_{r}h}B^{p}.
\end{align*}
Therefore, by assuming $\psi_{0}:=1$, we have that, for $k=1, 2, \ldots, p$,
\vspace{-2mm}
\begin{align}\label{The2}
\psi_{k-1}=e^{(k-1)\eta_{r}h}-\sum_{j=1}^{k-1}\varphi_{j}e^{(k-1-j)\eta_{r}h}.
\end{align}
\vspace{-3mm}
It follows from (\ref{The1}) that
\begin{align}\label{The3}
\phi(B)Y_n=\sum_{r=1}^{p}\big(Z_{n}^{(r)}+\psi_{1}Z_{n-1}^{(r)}+\ldots+\psi_{p-1}Z_{n-p+1}^{(r)}\big).
\end{align}
By replacing the noise $Z_n^{(r)}:=\alpha_r\int_{(n-1)h}^{nh}e^{\eta_r(nh-u)}dS_u$ and (\ref{The2}) in (\ref{The3}) we have that
\vspace{-1mm}
\begin{align*}
\phi(B)Y_n&=\sum_{k=1}^{p}\sum_{r=1}^{p}\big(e^{(k-1)\eta_{r}h}-\sum_{j=1}^{k-1}\varphi_{j}e^{(k-1-j)\eta_{r}h}\big)\alpha_r\int_{(n-k)h}^{(n-k+1)h}e^{\eta_r((n-k+1)h-u)}dS_u=:\sum_{k=1}^{p}V_{n-k+1}^{(k)}.
\end{align*}
\\
\\
\noindent
{\bf\large Appendix B}
\setcounter{equation}{0}
\begin{Proof}\label{B1}
\textnormal{: According to relation (\ref{X}) we have}
\end{Proof}
\vspace{-2mm}
\begin{align}\nonumber
{\bf X}^{*}_n&=\int_{-\infty}^{nh}e^{A(nh-u)}{\bf e}dS_u-E({\bf X}_n)=\int_{-\infty}^{(n-1)h}e^{A(nh-u)}{\bf e}dS_u
+\int_{(n-1)h}^{nh}e^{A(nh-u)}{\bf e}dS_u-E({\bf X}_n)\\ \nonumber
&=e^{Ah}\Big[\int_{-\infty}^{(n-1)h}e^{A((n-1)h-u)}{\bf e}dS_u\pm E({\bf X}_{n-1})\Big]+\int_{(n-1)h}^{nh}e^{A(nh-u)}{\bf e}dS_u-E({\bf X}_n)\\ \label{Xstar}
&=e^{Ah}{\bf X}^{*}_{n-1}+{\bf U}_{n},
\end{align}
where ${\bf U}_{n}:=\int_{(n-1)h}^{nh}e^{A(nh-u)}{\bf e}dS_u+e^{Ah}E({\bf X}_{n-1})-E({\bf X}_n)$.
Since $E({\bf X}^{*}_{n})={\bf 0}$ and $E({\bf X}^{*}_{n-1})={\bf 0}$, one can easily check from (\ref{Xstar}) that $E({\bf U}_n)={\bf 0}$. Furthermore, the covariance matrix of ${\bf U}_n$ is
\begin{align*}
Q_{n}:=\text{cov}({\bf U}_n,{\bf U}'_n)=\text{cov}\Big(\int_{(n-1)h}^{nh}e^{A(nh-u)}{\bf e}dS_u, \int_{(n-1)h}^{nh}{\bf e}'e^{A'(nh-u)}dS_u\Big).
\end{align*}
Now we find the covariance matrix ${\bf U}_{n}$ in two cases. First when $nh$ and $(n-1)h$ belong to one subinterval as $B_j$, so
\begin{align}\nonumber
Q_{n}&=E\Big[\Big(\int_{(n-1)h}^{nh}e^{A(nh-u)}{\bf e}dS_u-\int_{(n-1)h}^{nh}e^{A(nh-u)}{\bf e}E(dS_u)\Big)\\ \nonumber
&\quad\times\Big(\int_{(n-1)h}^{nh}{\bf e'}e^{A'(nh-u)}dS_u-\int_{(n-1)h}^{nh}{\bf e'}e^{A'(nh-u)}E(dS_u)\Big)\Big]\\ \nonumber
&=E\Big[\Big(\int_{(n-1)h}^{nh}e^{A(nh-u)}{\bf e}\big(dS_u-E(dS_u)\big)\Big)\times\Big(\int_{(n-1)h}^{nh}{\bf e'}e^{A'(nh-u)}\big(dS_u-E(dS_u)\big)\Big)\Big]\\ \label{Qn}
&=\int_{(n-1)h}^{nh}e^{A(nh-u)}{\bf e}{\bf e'}e^{A'(nh-u)}E\big(dS_u-E(dS_u)\big)^2=\int_{(n-1)h}^{nh}e^{A(nh-u)}{\bf e}
{\bf e'}e^{A'(nh-u)}{\text var}(dS_u).
\end{align}
Using (\ref{sl}) and the definition of $S_t$ in (\ref{LI}), the variance of the increment $dS_u$ for $u, u+du\in B_j$ is
\vspace{-1mm}
\begin{align}\nonumber
\text{var}(dS_u)&=\text{var}(S_{u+du}-S_{u})=\text{var}(S_{u+du}^{(1)}-S_{u}^{(1)})=\text{var}(S_{u+du}^{(1)})+\text{var}(S_{u}^{(1)})-2\text{cov}(S_{u+du}^{(1)},S_{u}^{(1)})\\ \label{varS}
&=\text{var}(S_{u+du}^{(1)})-\text{var}(S_{u}^{(1)})=\beta(\Lambda_{u+du}-\Lambda_{u})=\beta\frac{\lambda_j}{|B_{j}|}du,
\end{align}
and by changing the variable $nh-u$ to $u$ we have the relation (\ref{mat}).
Second, we consider $(n-1)h$ in $B_j$ and $nh$ in $B_{j+1}$, so
\vspace{-2mm}
\begin{align*}
Q_{n}=\text{cov}({\bf U}_n,{\bf U}'_n)=\text{cov}&\Big(\int_{(n-1)h}^{s_j}e^{A(nh-u)}{\bf e}dS_u+\int_{s_j}^{nh}e^{A(nh-u)}{\bf e}dS_u,\\
&\quad\int_{(n-1)h}^{s_j}{\bf e'}e^{A'(nh-u)}dS_u+\int_{s_j}^{nh}{\bf e'}e^{A'(nh-u)}dS_u\Big).
\end{align*}
Because of the independency of the increments $\text{cov}\Big(\int_{(n-1)h}^{s_{j}}e^{A(nh-u)}{\bf e}dS_u, \int_{s_{j}}^{nh}{\bf e'}e^{A'(nh-u)}dS_u\Big)=\text{cov}\Big(\int_{s_j}^{nh}e^{A(nh-u)}{\bf e}dS_u, \int_{(n-1)h}^{s_j}{\bf e'}e^{A'(nh-u)}dS_u\Big)=0$. Therefore,
\begin{align*}
Q_{n}=\text{cov}\big(\int_{(n-1)h}^{s_j}e^{A(nh-u)}{\bf e}dS_u, \int_{(n-1)h}^{s_j}{\bf e'}e^{A'(nh-u)}dS_u\big)
+\text{cov}\big(\int_{s_j}^{nh}e^{A(nh-u)}{\bf e}dS_u, \int_{s_j}^{nh}{\bf e'}e^{A'(nh-u)}dS_u\big),
\end{align*}
that by the same method in first case
\vspace{-1mm}
\begin{align*}
Q_{n}=\beta\frac{\lambda_j}{|B_{j}|}\int_{(n-1)h}^{s_j}e^{A(nh-u)}{\bf e}{\bf e'}e^{A'(nh-u)}du
+\beta\frac{\lambda_{j+1}}{|B_{j+1}|}\int_{s_j}^{nh}e^{A(nh-u)}{\bf e}{\bf e'}e^{A'(nh-u)}du.
\end{align*}
by changing the variable $nh-u$ to $u$ we have that
\vspace{-1mm}
\begin{align*}
Q_{n}=\beta\frac{\lambda_j}{|B_{j}|}\int_{nh-s_{j}}^{h}e^{Au}{\bf e}{\bf e'}e^{A'u}du
+\beta\frac{\lambda_{j+1}}{|B_{j+1}|}\int_{0}^{nh-s_{j}}e^{Au}{\bf e}{\bf e'}e^{A'u}du,
\end{align*}
so we get to the result of relation (\ref{mat}).
\vspace{2mm}
\begin{Proof}\label{B2}
\textnormal{: Since $\mathbf{X}^{*}_{1}=\mathbf{X}_{1}-E(\mathbf{X}_{1})$ is the centered state vector $\mathbf{X}_{1}:=\mathbf{X}_{h}$ and $\widehat{\mathbf{X}}_{1}^{*}={\bf 0}$, it follows from (\ref{X}) that}
\end{Proof}
\vspace{-7mm}
\begin{align*}
\Omega_{1}=\text{cov}(\mathbf{X}^{*}_{1},\mathbf{X}^{*^{'}}_{1})=\text{cov}(\mathbf{X}_{1},\mathbf{X}_{1}')=\text{cov}\Big(\int_{-\infty}^{h}e^{A(h-u)}
\mathbf{e} dS_u,\int_{-\infty}^{h}\mathbf{e}'e^{A'(h-u)}dS_u\Big).
\end{align*}
By a similar method in (\ref{Qn}), we have $\Omega_{1}=\int_{-\infty}^{h}e^{A(h-u)}\mathbf{e}\mathbf{e}'e^{A'(h-u)}\text{var}(dS_u)$. Therefore,
\begin{align*}
\Omega_{1}&=\int_{-\infty}^{0}e^{A(h-u)}\mathbf{e}\mathbf{e}'e^{A'(h-u)}\text{var}(dS_u)+\int_{0}^{h}
e^{A(h-u)}\mathbf{e}\mathbf{e}'e^{A'(h-u)}\text{var}(dS_u)\\
&=\sum_{k=0}^{\infty}\sum_{j=1}^{r}\int_{-s_{_{j}}-kT}^{-s_{_{j-1}}-kT}e^{A(h-u)}\mathbf{e}\mathbf{e}'e^{A'(h-u)}\text{var}(dS_u)+\int_{0}^{h}
e^{A(h-u)}\mathbf{e}\mathbf{e}'e^{A'(h-u)}\text{var}(dS_u).
\end{align*}
Similar to the relation (\ref{varS}), it follows from (\ref{sl}) and (\ref{Lam}) that  the variance of the increment $dS_u$ for $u\in[-s_{_{j}}-kT,-s_{_{j-1}}-kT)$ is
\vspace{-3mm}
\begin{align}\label{varS2}
\text{var}(dS_{u})=\text{var}(S_{u+du}-S_{u})=\text{var}(S_{-u}^{(2)}-S_{-u-du}^{(2)})=\beta(\Lambda_{-u}-\Lambda_{-u-du})
=\beta\frac{\lambda_j}{|B_{j}|}du.
\end{align}
The last equality follows from the fact that $-u\in(s_{_{j-1}}+kT, s_{_{j}}+kT]$. Let $h\in B_{1}$, then by (\ref{varS}) and (\ref{varS2}), we have
\begin{align*}
\Omega_{1}=\sum_{k=0}^{\infty}\sum_{j=1}^{r}\beta\frac{\lambda_{j}}{|B_{j}|}\int_{-s_{j}-kT}^{-s_{_{j-1}}-kT}e^{A(h-u)}\mathbf{e}\mathbf{e}'e^{A'(h-u)}du
+\beta\frac{\lambda_{1}}{|B_{1}|}\int_{0}^{h}e^{A(h-u)}\mathbf{e}\mathbf{e}'e^{A'(h-u)}du.
\end{align*}
By changing the variable $h-u$ to $u$, we have the relation $\Omega_{1}$.
\begin{Proof}\label{B3}
\textnormal{: According to relation $\mathbf{U}_{n}=\int_{(n-1)h}^{nh}e^{A(nh-u)}{\bf e}dL_u+e^{Ah}E({\bf X}_{n-1})-E({\bf X}_n)$, we have
\begin{align*}
Q_n=\text{cov}(\mathbf{U}_{n},\mathbf{U}'_{n})=\text{cov}\Big(\int_{(n-1)h}^{nh}e^{A(nh-u)}\mathbf{e}dL_{u},\int_{(n-1)h}^{nh}\mathbf{e}'e^{A'(nh-u)}dL_{u}\Big).
\end{align*}
By a similar method in (\ref{Qn}), it follows that
\begin{align*}
Q_n=\int_{(n-1)h}^{nh}e^{A(nh-u)}\mathbf{e}\mathbf{e}'e^{A'(nh-u)}\text{var}(dL_{u})=\text{var}(L_{1}) \int_{(n-1)h}^{nh}e^{A(nh-u)}\mathbf{e}\mathbf{e}'e^{A'(nh-u)}du.
\end{align*}
By changing the variable $nh-u$ to $u$, we have $Q_n=\text{var}(L_{1}) \int_{0}^{h}e^{Au}\mathbf{e}\mathbf{e}'e^{A'u}du$.}
\end{Proof}

\begin{Proof}\label{B4}
\textnormal{: Since $\mathbf{X}^{*}_{1}=\mathbf{X}_{1}-E(\mathbf{X}_{1})$ is the centered state vector $\mathbf{X}_{1}:=\mathbf{X}_{h}$ and $\widehat{\mathbf{X}}_{1}^{*}={\bf 0}$, we have that $\Omega_{1}=\text{cov}\big((\mathbf{X}_{1}^{*}-\widehat{\mathbf{X}}_{1}^{*})(\mathbf{X}_{1}^{*}-\widehat{\mathbf{X}}_{1}^{*})'\big)=\text{cov}(\mathbf{X}_{1},\mathbf{X}^{'}_{1}).$ So, it follows from Proposition 1 in \cite{bro5} that $\Omega_{1}=\text{var}(L_{1})\int_{0}^{\infty}e^{Au}\mathbf{e}\mathbf{e}^{'}e^{A^{'}u}du$.}
\end{Proof}

\end{document}